\documentclass{amsart}

\usepackage{graphicx}

\newtheorem{theorem}{Theorem}

\theoremstyle{definition}
\newtheorem{definition}{Definition}
\newtheorem{example}[theorem]{Example}

\theoremstyle{remark}

\numberwithin{equation}{section}

\begin{document}

\title{The Interlace Polynomial}

\author{Ada Morse}
\address{Department of Mathematics and Statistics, University of Vermont}
%\curraddr{Department of Mathematics and Statistics,
%Case Western Reserve University, Cleveland, Ohio 43403}
\email{ada.morse@uvm.edu}


%    General info
%\subjclass[2000]{Primary 54C40, 14E20; Secondary 46E25, 20C20}

\date{\today}


%\keywords{Differential geometry, algebraic geometry}

\begin{abstract}
In this paper, we survey results regarding the interlace polynomial of a graph, connections to such graph polynomials as the Martin and Tutte polynomials, and generalizations to the realms of isotropic systems and delta-matroids.
\end{abstract}

\maketitle

\section{Introduction}

The interlace polynomial of a graph arises in a number of settings. We begin with the simplest, that of a recursive method for counting Eulerian circuits in two-in two-out digraphs. The interlace polynomial of a simple graph is then obtained by generalizing the recursion used to solve this counting problem. We then discuss a closed form for the polynomial in terms of its adjacency matrix, the structure of which suggests definitions for analogous polynomials as well as a two-variable generalization. Another context in which the interlace polynomial arises is in isotropic systems, where it appears as a specialization of the Tutte-Martin polynomials, a connection we follow by way of the Martin polynomials of $4$-regular graphs. Lastly, we review generalizations of the polynomial to square matrices and delta-matroids.

In the context of counting Eulerian circuits in two-in two-out digraphs, the interlace polynomial arose by way of Arratia, Bollob\'{a}s, and Sorkin's work on DNA sequencing \cite{ABS2000}. In DNA sequencing by hybridization, the goal is to reconstruct a string of DNA knowing only information about its shorter substrings. The problem is to determine, from knowledge about the shorter substrings, whether a unique reconstruction exists.

More precisely, if $A = a_1a_2\cdots a_m$ is a sequence consisting of $m$ base pairs, the $l$-spectrum of $A$ is the multiset containing all $l$-tuples consisting of $l$ consecutive base pairs in $A$. Given knowledge of the $l$-spectrum, the goal is to determine the number $k_l(m)$ of sequences of base pairs of length $m$ having that $l$-spectrum.

In \cite{ABS2000}, the authors associate to a given $l$-spectrum its \emph{de Bruijn graph}: a two-in two-out digraph $D$ such that the Eulerian circuits of $D$ are in bijection with sequences of base pairs having that $l$-spectrum. The problem, then, is to count the number of Eulerian circuits of $D$. This approach led to the discovery of a recursive formula for computing the number of Eulerian circuits of $D$ based on an associated interlace graph. In \cite{ABS2004ip}, Arratia, Bollob\'{a}s, and Sorkin generalized this recursion to define the interlace polynomial of an arbitrary simple graph.

The Eulerian circuits and cycle decompositions of $4$-regular graphs have been an area of significant interest among graph theorists for many years, and approaches using graph polynomials have frequently proved fruitful \cite{M77,LV83, J90, EMS2002}. The Martin polynomial \cite{M77,LV83}, in particular, is closely related to the interlace polynomial as it counts, for any $k$, the number of $k$-component circuit partitions of a $4$-regular graph. 

This connection can be made explicit and, indeed, generalized. In a series of papers in the 1980s-1990s, Bouchet introduced the notion of an isotropic system to unify aspects of the study of $4$-regular graphs and binary matroids \cite{B87,B88,B91}, including a generalization of the Martin polynomials to this area \cite{B91}. Shortly after the discovery of the interlace polynomial, it was noticed that the interlace polynomial can be found as a specialization of the (restricted) Tutte-Martin polynomial of an isotropic system \cite{B05,AH2004}.

A connection between the interlace polynomial and the Tutte polynomial can be found by way of the Martin polynomial. However, this connection only captures the Tutte polynomial $t(G;x,y)$ for plane graphs when $x=y$, and so does not provide any strong link between the interlace polynomial and the many specializations of the Tutte polynomial, such as the chromatic polynomial. Slightly more of the Tutte polynomial can be captured by a generalization of the interlace polynomial of a matroid, however there is no known connection between the interlace polynomial and the full Tutte polynomial (or other deletion-contraction polynomial) in any context.

Many generalizations of the interlace polynomial have been obtained. In \cite{ABS2004tvip}, Arratia, Bollob\'{a}s, and Sorkin defined a two-variable interlace polynomial of which the single-variable polynomial is a specialization. In doing so, they discovered, concurrently with Aigner and van der Holst \cite{AH2004}, a closed form for the single-variable interlace polynomial in terms of its adjacency matrix. This closed form has a natural extension to arbitrary square matrices, and using  a delta-matroid associated to the adjacency matrix of a graph, Brijder and Hoogeboom obtained a generalization of the interlace polynomial to delta-matroids \cite{BH2014ipm}. In each case, the recursive definition of the interlace polynomial has also been generalized.

\section{The interlace polynomial of a graph}

We first begin by defining the interlace polynomial recursively by way of counting Eulerian circuits in two-in two-out digraphs, and then discuss a closed form, analogous polynomials, and a two variable generalization. We conclude with selected evaluations of the interlace polynomial.

\subsection{Preliminary definitions}

\begin{figure}
    \centering
    \includegraphics{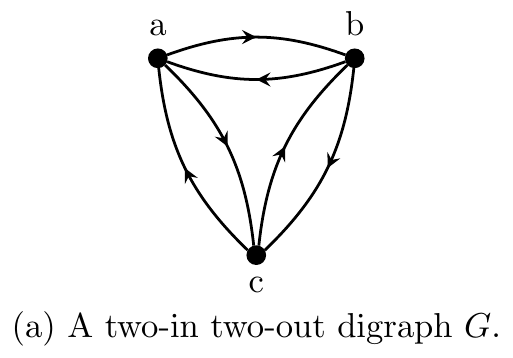}
    \includegraphics{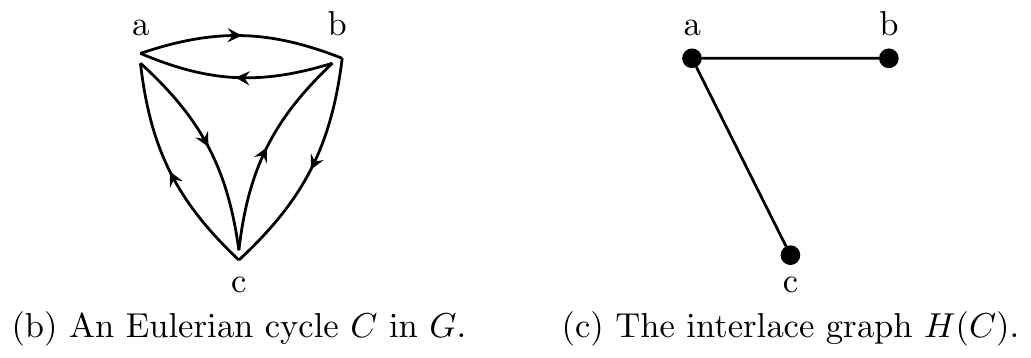}
    \includegraphics{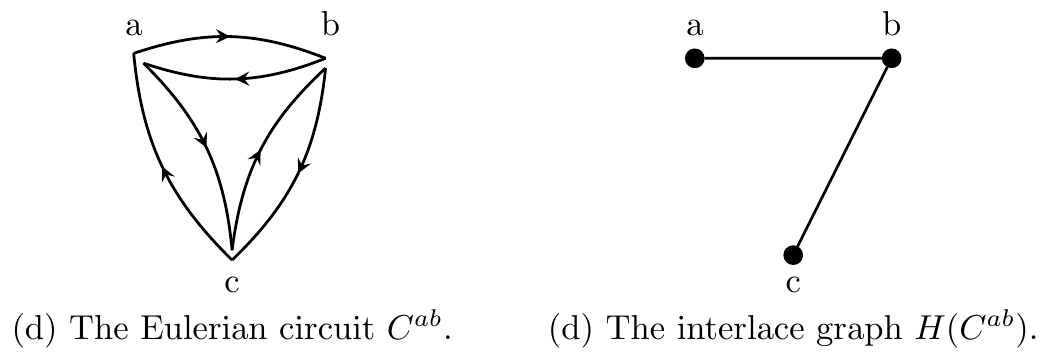}
    \caption{Transpositions of Eulerian circuits and the interlace graph.}
    \label{fig:trans}
\end{figure}

\begin{definition} \label{def:interlace-graph} A \emph{two-in two-out   digraph} is a $4$-regular digraph such that each vertex has both indegree and outdegree equal to two. Let $G$ be a two-in two-out digraph. An \emph{Eulerian circuit} of $G$ is a closed, directed walk of $G$ containing each edge exactly once. Given an Eulerian circuit $C$ of $G$, we say that vertices $a$ and $b$ are \emph{interlaced} if the cycle visits them in the order $a\ldots b \ldots a \ldots b \ldots$ and \emph{noninterlaced} otherwise. The \emph{interlace graph} or \emph{circle graph} $H(C)$ is the graph whose vertices are the vertices of $G$ with edges between pairs of vertices interlaced in $C$ (see Figure \ref{fig:trans} (b) and (c)).
\end{definition}

Interlace graphs have been extensively studied \cite{B87circle,B87circle2,B94,F97,RR76,R99}, and were characterized by Bouchet in \cite{B87circle2,B94}. A particular focus of the area, due to a problem of Gauss, has been characterizations of the interlace graphs arising from Eulerian cycles in plane $4$-regular graphs \cite{F97,RR76,R99}.

There is a natural operation defined on Eulerian circuits of two-in two-out digraphs in terms of this interlace relation.

\begin{definition} \label{def:transposition} At each vertex of a two-in two-out digraph $G$, there are two possible (orientation consistent) pairings of in-edges and out-edges. For a pair of vertices $a$ and $b$ interlaced in an Eulerian circuit $C$ of $G$, define the \emph{transposition} $C^{ab}$ to be the Eulerian circuit obtained by switching the pairing of edges at $a$ and $b$, see Figure \ref{fig:trans} (b) and (d).
\end{definition}

The Eulerian circuits of $G$ form a single orbit under the action of transposition, the proof of which can be found in \cite{ABS2004ip} but was known previously in more general form in \cite{P95,U92}. Observation of the effect on the interlace relation by performing the above operation to an Eulerian circuit leads to a corresponding definition for interlace graphs, presented here for graphs in general.

\begin{figure}
    \centering
    \includegraphics{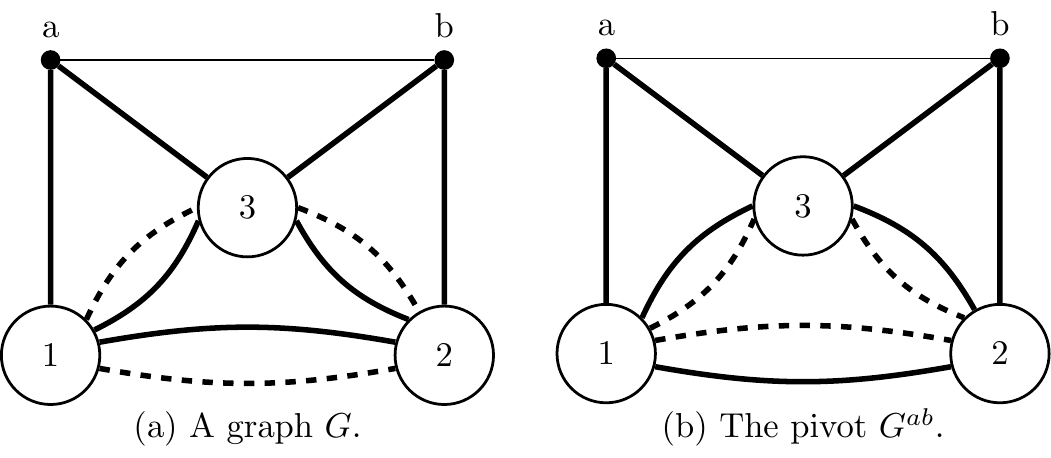}
    \caption{On the left, a graph $G$ with edge $ab$ and vertices partitioned as in Definition \ref{def:pivot} (the parts of $G$ unaffected by pivoting are not shown.) On the right, $G^{ab}$ is obtained by toggling edges/nonedges among the sets of vertices labeled $1, 2,$ and $3$.}
    \label{fig:pivot}
\end{figure}

\begin{definition} \label{def:pivot} Let $G$ be any graph. Let $v \in V(G)$. For any pair of vertices $a,b \in V(G)$, partition the remaining vertices of $G$ into the following sets: (1) vertices adjacent to $a$ and not $b$; (2) vertices adjacent to $b$ and not $a$; (3) vertices adjacent to both $a$ and $b$; and (4) vertices adjacent to neither $a$ nor $b$. Define the \emph{pivot} $G^{ab}$ to be the graph obtained by inserting all possible edges between the first three of these sets, and deleting those that were already present in $G$ (see Fig. \ref{fig:pivot}). Denote by $G_{ab}$ the graph $G$ with the labels of the vertices $a$ and $b$ swapped.
\end{definition}

While the definition of pivot above is due to Arratia, Bollob\'{a}s, and Sorkin \cite{ABS2000}, the idea of the pivot appeared in the earlier work of Kotzig \cite{K68} on local complementations and the graph $(G^{ab})_{ab}$ is defined by Bouchet in \cite{B88} as the \emph{complementation of $G$ along the edge $ab$}. The precise connection to both is as follows.

\begin{definition} \label{def:graph} Let $G$ be a graph. For $v \in V(G)$, we denote the open neighborhood of $v$ by $N(v)$. We say $v$ is a \emph{looped vertex} if $v$ has a loop. Note that $v \not \in N(v)$ even if $v$ is a looped vertex. We define the \emph{local complement} $G*v$ to be the graph obtained from $G$ by interchanging edges and non-edges in $N(v)$. By convention, we read graph operations left-to-right, so $G*v*w*v = ((G*v)*w)*v$.
\end{definition}

\begin{theorem} \label{thm:localcomp} \cite{B01,B88} Let $G$ be a graph. If $ab$ is an edge in $G$ with neither $a$ nor $b$ a looped vertex, then $(G^{ab})_{ab} = G*a*b*a$.
\end{theorem}

In the case of interlace graphs, the pivot operation captures the behavior of a transposition of an Eulerian circuit in the following sense.

\begin{theorem} \label{prop:pivottransposition} \cite{ABS2004ip} For an Eulerian circuit $C$ of a two-in two-out digraph $G$, we have $(H(C))^{ab} = (H(C^{ab}))_{ab}$.
\end{theorem}

We can now define the interlace polynomial of a graph. Arratia, Bollob\'{a}s, and Sorkin proved in \cite{ABS2004ip} that the recurrence below does not depend on the order of edges chosen, i.e. the polynomial is well-defined. 

\begin{definition} \emph{(The interlace polynomial)} \label{thm:ipolynomial} \cite{ABS2004ip} Let $G$ be a simple graph. The \emph{the interlace polynomial of $G$}, denoted $q_N(G;x)$, is defined by
\begin{displaymath}
   q_N(G;x) = \left\{
     \begin{array}{lr}
       q_N(G \setminus a;x) + q_N(G^{ab} \setminus b;x), & ab \in E(G)\\
       x^n, & G \cong E_n
     \end{array}
   \right.
\end{displaymath}
where $E_n$ is the graph on $n$ vertices with no edges.
\end{definition}

Note that while the recurrence above is presented in its original form, in generalizations of the interlace polynomial the label-switching operation $G_{ab}$ (see Definition \ref{def:pivot}) occurs as part of the generalized pivot operation. In the case of the recurrence above, this can be obtained using local complementation in place of the pivot operation (see Theorem \ref{thm:localcomp}). Under that convention, the recurrence above becomes $q_N(G;x) = q_N(G \setminus a;x) + q_N(G*a*b*a \setminus a;x)$, which aligns with the form of the recurrence used in subsequent sections. In addition, the interlace polynomial was originally denoted $q(G)$. We follow \cite{ABS2004tvip} in reserving that notation for the two-variable generalization.

Definition \ref{thm:ipolynomial} is stated for simple graphs. It can, however, be extended to the case of looped graphs (i.e. graphs with loops, but without multiple edges.) In this case, the recurrence above only holds for edges where neither endpoint has a loop, and an additional recurrence is required to handle looped vertices. For precise details, see Theorem \ref{thm:tviprecurrence} below on the two-variable interlace polynomial.

Aigner and van der Holst discovered in \cite{AH2004} a state-sum formulation for the interlace polynomial in terms of the adjacency matrix. 

\begin{definition} \label{def:adjacencymatrix}  Let $A(G)$ be the adjacency matrix of a graph $G$. For $T \subseteq V(G)$, we denote by $G[T]$ the subgraph of $G$ induced by $T$. We denote by $n(G[T])$ and $r(G[T])$ the nullity and rank respectively of $A(G[T])$. By convention, $n(G[\emptyset]) = 0$. For $S \subseteq V(G)$, we define the \emph{loop complement} of $G$ with respect to $S$, denoted $G+S$, to be the graph obtained by adding loops to unlooped vertices in $S$ and removing loops from looped vertices of $S$.
\end{definition}

\begin{theorem} \label{thm:vertexnullity} \cite{AH2004} Let $G$ be a simple graph. Then
    \begin{equation*}
        q_N(G;x) = \sum_{T \subseteq V(G)} (x-1)^{n(G[T])}.
    \end{equation*}
\end{theorem}

Due to this formula, the interlace polynomial is sometimes referred to as the vertex-nullity polynomial. A related polynomial is the vertex-rank polynomial, obtained by replacing $n(G[T])$ with $r(G[T])$ in the above expression (see the discussion of the two-variable interlace polynomial below.)

Aigner and van der Holst \cite{AH2004} as well as Bouchet \cite{B05} also defined and studied the following related polynomial.

\begin{definition}\label{def:globalinterlace} \cite{AH2004,B05} Let $G$ be a graph with vertex $a$. Define the polynomial $Q(G;x)$ by the following recursion:
\begin{enumerate}
    \item if $G$ contains an edge $ab$ then
        \begin{equation*}
            Q(G;x) = Q(G \setminus a;x) + Q(G*a \setminus a;x) + Q(G^{ab}\setminus b;x), and
        \end{equation*}
    \item $Q(E_n;x) = x^n$.
\end{enumerate}
\end{definition}

As with the interlace polynomial, there is also an adjacency matrix version of the $Q$ polynomial. Note that the expression below is not the original form given in \cite{AH2004}, but can be recovered from the description in \cite{AH2004} (see also Traldi's work on labeled interlace polynomials \cite{T13}).

\begin{theorem} \cite{AH2004} \label{thm:globalstatesum} For a graph $G$ we have
    \begin{equation*}
        Q(G;x) = \sum_{T \subseteq V(G)} \sum_{S \subseteq T} (x-2)^{n((G+S)[T])}.
    \end{equation*}
\end{theorem}

In \cite{ABS2004tvip}, Arratia, Bollob\'{a}s, and Sorkin developed a two-variable extension of the interlace polynomial. Note that while we will focus on the two-variable polynomial below, another multivariable generalization was studied by Courcelle in \cite{C2008} and Traldi studied a labeled multivariable interlace polynomial in \cite{T13}.

\begin{definition} \label{def:tvip} Let $G$ be a graph. Then the two-variable interlace polynomial is
    \begin{equation*}
    q(G;x,y) = \sum_{T \subseteq V(G)} (x-1)^{r(G[T])}(y-1)^{n(G[T])}
\end{equation*}
\end{definition}

This is indeed an extension of the single-variable polynomial: setting $x=2$ in the above equation yields precisely the formula of Theorem \ref{thm:vertexnullity}. Setting $y=2$ instead yields a related graph polynomial (the \emph{vertex-rank polynomial}), studied in \cite{ABS2004tvip}. The two-variable polynomial also satisfies a recurrence generalizing that satisfied by the single-variable polynomial. Indeed, on simple graphs, setting $x=2$ in the recurrence below recovers the original recurrence of the single-variable interlace polynomial. Thus, on looped graphs, setting $x=2$ provides an extension of the single-variable polynomial to graphs with loops.

\begin{theorem}\cite{ABS2004tvip}\label{thm:tviprecurrence} The two-variable interlace polynomial satisfies the following recurrence:
\begin{enumerate}
    \item if $ab$ is an edge of $G$ where neither $a$ nor $b$ has a loop, then
        \begin{equation*}
            q(G;x,y) = q(G\setminus a;x,y) + q(G^{ab}\setminus b;x,y) + ((x-1)^2 - 1)q(G^{ab} \setminus a \setminus b;x,y),
        \end{equation*}
    \item if $a$ is a looped vertex of $G$ then
        \begin{equation*}
            q(G) = q(G \setminus a;x,y) + (x-1)q(G * a \setminus a;x,y), and
        \end{equation*}
    \item $q(E_n;x,y) = y^n$.
\end{enumerate}
\end{theorem}

Lastly, note that while the single-variable interlace polynomial can be viewed as a specialization of the restricted Tutte-Martin polynomial of an isotropic system (see Section \ref{sec:isotropic} below), no such generalization is known for the two-variable version.

\subsection{Evaluations of the interlace polynomial}

The interlace polynomial of a graph has been found to encode structural information as well as graph invariants. These include Eulerian circuits, perfect matchings, independence number, component number, and more. The evaluations below of $q_N(G)$ at $1, -1, 3$ and $2$ extend to graphs with loops while the evaluation at $0$ does not \cite{BH2014ipm}. The evaluation at $-1$ was conjectured in \cite{ABS2000}. The proofs of these evaluations specifically for the interlace polynomial can be found in the papers cited below, but we note that many can be recovered from evaluations of the Tutte-Martin polynomials derived in \cite{B91}. Item (2) is the solution to the counting problem that motivated the development of the polynomial.

\begin{theorem} \label{thm:qeval} Let $G$ be a graph, possibly with loops but without multiple edges. Let $n = |V(G)|$.
\begin{enumerate}
    \item \cite{AH2004} $q_N(G;1)$ is the number of induced subgraphs of $G$ with an odd number of perfect matchings (including the empty set.)
    \item \cite{ABS2004ip} If $H(C)$ is the interlace graph of an Eulerian circuit of a two-in two-out digraph $D$, then $q_N(H(C);1)$ is the number of Eulerian circuits in $D$.
    \item \cite{ABS2004ip} $q_N(G;2) = 2^n$.
    \item \cite{AH2004,BBCP2002,BH2014ipm} $q_N(G;-1) = (-1)^n(-2)^{n(G+V(G))}$.
    \item  \cite{AH2004} If $G$ is simple then $q_N(G;0) = 0$ if $n \geq 1$.
    \item \cite{AH2004} $q_N(G;3) = kq_N(G;-1)$ for some odd integer $k$.
\end{enumerate}
\end{theorem}

\begin{theorem} \label{thm:globalqevals} \cite{AH2004}  Let $G$ be a simple graph with $n=|V(G)|$.
\begin{enumerate}
    \item $Q(G;0) = 0$ if $n\geq 1$.
    \item $Q(G;3) = 3^n$.
    \item $Q(G;4) = 2^n e$ where $e$ is the number of induced Eulerian subgraphs of $G$.
    \item For each $T \subseteq V(G)$, we associate with $T$ general induced subgraphs, which are subgraphs obtained from the subgraph induced by $T$ by adding loops at any of the vertices of $T$. We allow perfect matchings of a general induced subgraph to include loops. Then $Q(G)(2)$ is the number of general induced subgraphs with an odd number of general perfect matchings.
\end{enumerate}
\end{theorem}

The following results describe the structure of the interlace polynomial.

\begin{theorem} Let $G$ be a simple graph with $n=|V(G)|$.
\begin{enumerate}
    \item \cite{AH2004} Let $[G]$ denote the set of all graphs obtainable from $G$ by a sequence of pivots. Then $\deg q_N(G;x) = \max_{H \in [G]} \alpha(H)$ where $\alpha(H)$ is the independence number of $H$.
    \item \cite{ABS2004ip} The least power of $x$ appearing in $q_N(G;x)$ is the number of components of $G$.
    \item \cite{EMS07} If $n \geq 1$ then $q_N(G;x)$ has no constant term.
    \item \cite{EMS07} If $n>1$ then writing $q_N(G;x) = \sum a_i x^i$ and $q(G;x,y) = \sum a_{ij} x^i y^j$ yields $a_1 = a_{01} = -a_{10}$.
    \item \cite{EMS07} If $n>1$ then writing $q(G;x,y)$ and $q_N(G;x)$ as above we have $a_1 = \sum_{i \geq 0} a_{i1} 2^i$ and $\sum_{i \geq 1} a_{i1}2^i = 0$.
\end{enumerate}
\end{theorem}

The common value $a_1 = a_{01} = -a_{10}$ in item (4) above is defined and studied as a graph invariant in \cite{EMS07}.

\section{Connections to other polynomials}

Graph polynomials have been used extensively in the study of $4$-regular graphs and their circuit decompositions, and there are many connections between the interlace polynomial and other graph polynomials arising in that context.

\subsection{The Martin and Tutte polynomials}

The Martin polynomial was defined by Martin in \cite{M77} to study circuit partitions of $4$-regular graphs. Given that the interlace polynomial can be used to count the number of Eulerian circuits of a two-in two-out digraph, it is not surprising that the polynomials should have some connection. Before defining the Martin polynomial(s), we establish the following notions for $4$-regular graphs.

\begin{definition} Let $G$ be a $4$-regular graph. A \emph{circuit partition} of $G$ is a decomposition of $G$ into edge-disjoint circuits. A \emph{transition} (or \emph{state}) at a vertex $v$ of $G$ is a choice of one of the three possible pairings of edges incident with $v$. If $G$ is a two-in two-out digraph, we require that transitions follow the orientation of $G$ by pairing incoming edges with outgoing edges. A \emph{transition system} (or \emph{graph state}) $T$ of $G$ consists of a choice of transition at each vertex of $G$. Any transition system $T$ of $G$ induces a circuit partition of $G$ and vice versa. Denote by $|T|$ the number of circuits in the transition system $T$. An \emph{Eulerian system} is a choice of Eulerian circuit for each component of $G$.
\end{definition}

\begin{figure}
    \centering
    \includegraphics{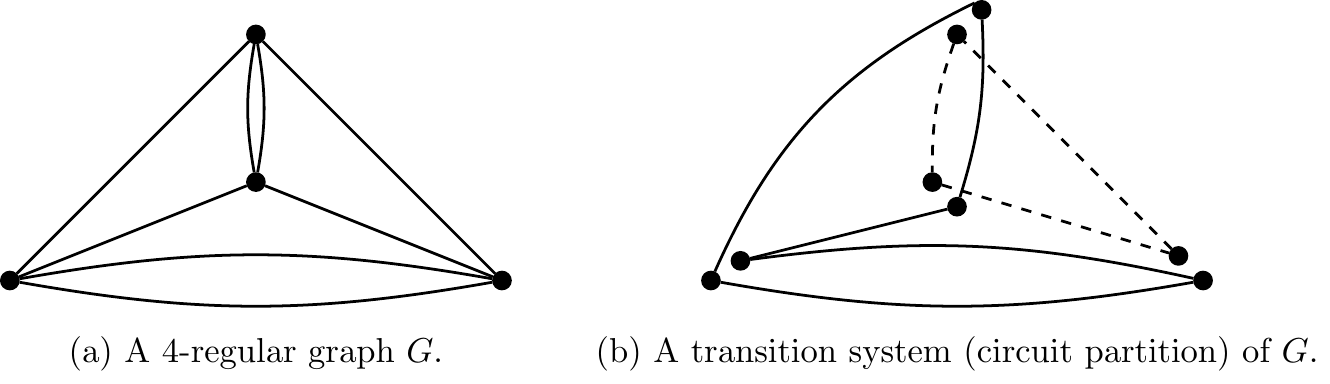}
    \caption{A transition system of a $4$-regular graph.}
    \label{fig:transition-system}
\end{figure}

We note that the Martin polynomials given below were originally defined recursively by Martin in \cite{M77}, with the closed forms due to Las Vergnas \cite{LV83}. We follow closely the notation of \cite{B91} to make the connection to Tutte-Martin polynomials in the next section most explicit.

\begin{definition} \label{def:martindigraph} Let $G$ be a $4$-regular graph. The \emph{Martin polynomial of $G$} is
\begin{equation*}
    M(G;x) = \sum (x-2)^{|T| - k(G)}
\end{equation*}
where the sum is over all transition systems $T$. For a two-in two-out digraph $G$ we define
\begin{equation*}
    m(G;x) = \sum (x-1)^{|T| - k(G)}
\end{equation*}
where the sum is over all transition systems $T$.
\end{definition}

The connection of these polynomials to the interlace polynomials $q_N$ and $Q$ of graphs can be seen as follows. Let $G$ be a $4$-regular graph and let $C$ be an Eulerian system of $G$. Let $H(C)$ be the interlace graph of $C$. Let $P$ be any circuit partition of $G$. At each vertex, the transition in $P$ is either contained in $C$, is consistently oriented by $C$ but not contained in $C$, or is not consistently oriented by $C$. Let $W$ be the set of vertices at which $P$ agrees with $C$, $Y$ the set at which $P$ disagrees but follows the orientation induced by $C$, and $Z$ the set of vertices at which $P$ disagrees with the orientation induced by $C$. Then Traldi has proven the following result in \cite{T11}, building on work of Cohn and Lempel \cite{CL72}.

\begin{theorem} \cite{T11} \label{thm:cohnlempel}Under the conditions of the previous remarks, we have 
    \begin{equation} \label{eq:cohnlempel}
        |P| - k(G) = n((H(C)+Z)[Y \cup Z])
    \end{equation}
\end{theorem}

This connection between circuit partitions and nullities yields the following equality between the Martin  and interlace polynomial, a result initially observed in \cite{ABS2004ip} and proved in \cite{EMS07,T11}.

\begin{theorem} \label{thm:martinq} Let $G$ be a two-in two-out digraph and let $C$ be an Eulerian decomposition of $G$. Then $m(G;x) = q(H(C);x)$.
\end{theorem}

Equation \ref{eq:cohnlempel} can also be used to obtain the following theorem.

\begin{theorem} Let $G$ be a $4$-regular graph. Let $C$ be an Eulerian system of $G$. Then $M(G;x) = Q(H(C);x)$.
\end{theorem}

Theorem \ref{thm:cohnlempel} can also be used to obtain a connection between the interlace polynomial and the Tutte polynomial. We recall here the recursive definition of the Tutte polynomial, and refer the reader to e.g. \cite{EMM11,BO92} for surveys.

\begin{definition} Let $G$ be a graph. The \emph{Tutte polynomial} of $G$ is the polynomial $t(G;x,y)$ obtained from the following recurrence:
\begin{enumerate}
    \item $t(G;x,y) = t(G\setminus e;x,y) + t(G / e; x,y)$ if $e$ is an edge that is neither a bridge nor a loop, and
    \item $t(G;x,y) = x^i y^j$ if no such edge exists and $G$ has $i$ bridges and $j$ loops.
\end{enumerate}
\end{definition}

Martin discovered a connection between the Martin polynomial and the Tutte polynomial in the case of plane graphs, which then extends, by results above, to the interlace polynomial. We first require the definition of the medial graph of a plane graph.

\begin{figure}
    \centering
    \includegraphics{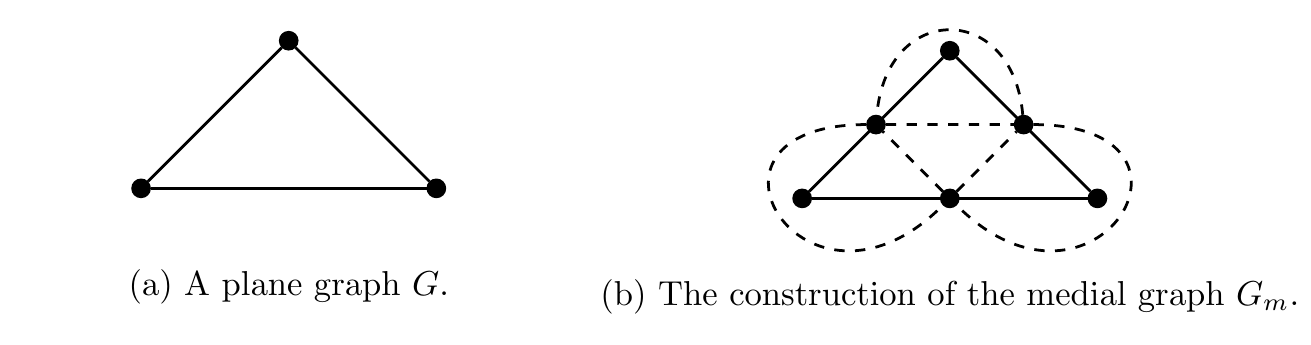}
    \includegraphics{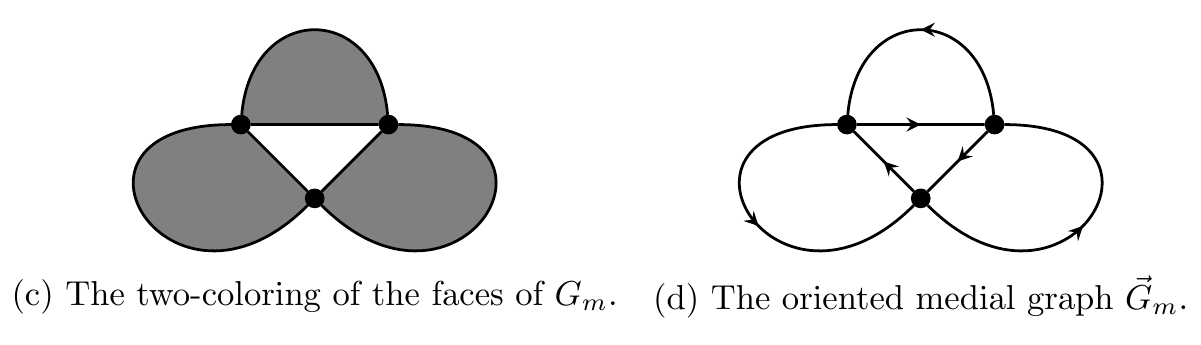}
    \caption{The construction of the oriented medial graph $\vec{G}_m$ from a plane graph $G$.}
    \label{fig:medial}
\end{figure}

\begin{definition} \label{def:medial} Let $G$ be a plane graph. The medial graph $G_m$ of $G$ is obtained by placing vertices on each of the edges of $G$, and connecting these vertices with edges by following the face-boundary walks (see Fig. \ref{fig:medial}). Color the faces of $G_m$ containing a vertex of $G$ black, and color the remaining faces white. This properly two-colors the faces of $G$. Orient edges counterclockwise around the black faces to obtain the oriented medial graph $\vec{G}_m$. Note that $\vec{G}_m$ is a two-in two-out digraph. \end{definition}

With this construction, we have the following two theorems relating the Martin, Tutte, and interlace polynomials.

\begin{theorem} \cite{M78} Let $G$ be a plane graph with oriented medial graph $\vec{G}_m$. Then 
    \begin{equation*}
        t(G;x,x) = m(\vec{G}_m;x).
    \end{equation*}
\end{theorem}

\begin{theorem} \label{thm:tutteplane} \cite{EMS07} Let $G$ be a plane graph with oriented medial graph $\vec{G}_m$. Let $C$ be an Eulerian circuit in $\vec{G}_m$ with interlace graph $H(C)$. Then
    \begin{equation*}
        t(G;x,x) = q_N(H(C);x).
    \end{equation*}
\end{theorem}

\subsection{Isotropic systems and the Tutte-Martin polynomials}
\label{sec:isotropic}

Isotropic systems were introduced and studied by Bouchet in a number of papers to unify the study of binary matroids and transition systems of $4$-regular graphs \cite{B87,B88,B91}. In particular, he introduced in \cite{B91} the Tutte-Martin polynomials of isotropic systems, of which the one-variable interlace polynomial of a graph is a specialization. We follow here the notation and approach of \cite{B05}.

\begin{definition} \label{def:isotropic} Let $K = \{0,x,y,z\}$ be the Klein $4$-group under addition, considered as a vector space of dimension $2$ over $GF(2)$. Let $K' = K \setminus 0$. Let $\langle \cdot, \cdot \rangle$ be the bilinear form on $K$ given by $\langle a,b \rangle = 1$ if neither $a$ nor $b$ is zero and $a \neq b$ and $\langle a,b \rangle = 0$ otherwise. For any finite set $V$, denote by $K^V$ the set of $V$-tuples with entries from $K$ considered as a vector space over GF$(2)$. Define $(K')^V$ similarly. Extend the bilinear form on $K$ to $K^V$ by defining $\langle X,Y \rangle = \sum_{v \in V} \langle X_v, Y_v \rangle$ (where e.g. $X_v$ is the entry in the $v$-labelled coordinate of $X$). We define an \emph{isotropic system} to be a pair $(V,L)$ where $V$ is a finite set, $L$ is a subspace of $K^V$ of dimension $|V|$, and $\langle X,Y \rangle = 0$ for all $X,Y \in L$.  For any $X \in K^V$, define $\widehat{X} = \{Y \in K^V : Y_v \in \{0,X_v\} \text{ for all } v \in V\}$. Note that $\widehat{X}$ is always an isotropic system.
\end{definition}

\begin{definition} \label{def:tutte-martin} Let $S = (V,L)$ be an isotropic system. Let $C \in (K')^V$. The \emph{restricted Tutte-Martin polynomial} of $S$ with respect to $C$ is given by
    \begin{equation*}
        tm(S,C;x) = \sum_X (x-1)^{\dim(L \cap \widehat{X})}
    \end{equation*}
where the sum is taken over all $X \in (K')^V$ such that $X_v \neq C_v$ for all $v$.
\end{definition}

\begin{figure}
    \centering
    \includegraphics{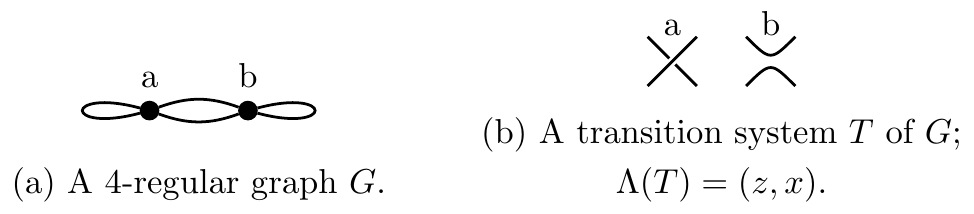}
    \includegraphics{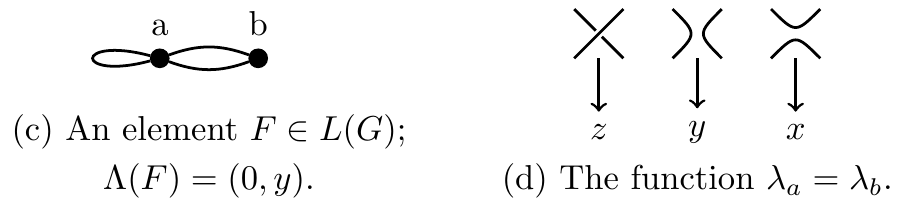}
    \caption{Construction of the isotropic system associated to a $4$-regular graph $G$ and labelling function $\lambda_v$. In this case, the image $L$ of $L(G)$ under $\Lambda$ is the isotropic system $\{ (0,0), (0,y), (y,0), (y,y)\}$.}
    \label{fig:isotropic}
\end{figure}

The connection between isotropic systems and $4$-regular graphs can be seen as follows (see \cite{B91}). Let $G = (V,E)$ be a $4$-regular graph. The \emph{cycle space} $L(G)$ is the collection of all edge-sets of $G$ inducing subgraphs having even degree at each vertex. For each vertex $v \in V$, let $\lambda_v$ be a bijection labelling the three transitions at $v$ with distinct values from $K'$ (see Figure \ref{fig:isotropic} (d)). This labelling induces a bijection $\Lambda$ from transition systems of $G$ to $(K')^V$ by defining $\Lambda(T)_v = \lambda_v(t)$ where $t$ is the transition of $T$ at the vertex $v$ (see Figure \ref{fig:isotropic} (b)). The labelling $\lambda$ can also be used to define a linear map from $L(G)$ to $K^V$ as follows. Given $F \in L(G)$ and $v \in V$, $F$ contains either no edges at $v$, four edges at $v$, or two edges at $v$. In the first and second case, define $\Lambda(F)_v = 0$. In the final case, the pairing of those two edges induces a transition $t$ at $v$, and we define $\Lambda(F)_v = \lambda_v(t)$ (see Figure \ref{fig:isotropic} (c)). The image $L$ of $L(G)$ under $\Lambda$ is a subspace of $K^V$. Bouchet has proven \cite{B87} that $S = (L,V)$ is an isotropic system such that $|T| - k(G) = \dim(L \cap \widehat{\Lambda(T)})$ for any transition system $T$. This yields the following connection to the Martin polynomial.

\begin{theorem} \cite{B91} \label{thm:martintm} Let $G$ be a two-in two-out digraph with transition system $T$. Then
    \begin{equation*}
        tm(S,\Lambda(T);x) = m(G;x).
    \end{equation*}
\end{theorem}

While the above provides an immediate connection to the interlace polynomial on interlace graphs via Theorem \ref{thm:martinq}, we can in fact recapture the interlace polynomial for any graph through a more general construction of an isotropic system associated to a graph.

\begin{definition} \label{ex:graphic-presentation} Let $G$ be a simple graph with vertex set $V$. Let $N(v)$ denote the neighborhood of $v \in V$. The powerset of $V$ forms a vector space over $GF(2)$ with addition corresponding to symmetric difference of sets. For $P \subseteq V$, define $N(P) = \sum_{v \in P} N(v)$. For $X \in K^V$ and $P \subseteq V$, denote by $X|P$ the vector given by $(X|P)_v = X_v$ for $v \in P$ and $(X|P)_v = 0$ for $v \not \in P$. Let $A,B \in (K')^V$ with $A_v \neq B_v$ for all $v \in V$. Define $L = \{ A|P + B|(N(P)) : P \subseteq V\}$. Then $S=(V,L)$ is an isotropic system, for which the triple $(G,A,B)$ is called a \emph{graphic presentation}.
\end{definition}

Aspects of the above definition can be seen as generalizations of the $4$-regular case. For example, in the $4$-regular case, when $S = (L,V)$ is associated as in the discussion above to the pair $(G,\{\lambda_v:v \in V(G)\})$, Eulerian circuits of $G$ correspond to vectors $X \in (K')^V$ such that $\dim(L \cap \widehat{X}) = 0$. On the other hand, when $(G,A,B)$ is a graphic presentation for $S = (L,V)$, the function $A$ satisfies $\dim(L \cap \widehat{A}) = 0$, so $A$ in some sense generalizes Eulerian circuits to this case. For a detailed exposition, see \cite{B88}.

\begin{theorem} \label{thm:tutte-martin} \cite{B05} Let $G$ be a simple graph and let $S$ be the isotropic system associated to the graphic presentation $(G,A,B)$. Then
    \begin{equation*}
        tm(S,A+B;x) = q_N(G;x).
    \end{equation*}
\end{theorem}

The $Q$ polynomial of Definition \ref{def:globalinterlace} is also a specialization of a polynomial of isotropic systems. Note that once again a natural connection arises from the construction preceding Theorem \ref{thm:martintm}. However, in this case we can, as for the interlace polynomial, recover $Q$ entirely from graphic presentations.

\begin{definition} \label{def:globaltuttemartin} \cite{B05} The \emph{global Tutte-Martin polynomial} of an isotropic system $S = (L,V)$ is
    \begin{equation*}
        TM(S;x) = \sum_{X \in (K')^V} (x-2)^{\dim(L \cap \widehat{X})}.
    \end{equation*}
\end{definition}

\begin{theorem} \label{thm:globalrelation} \cite{AH2004,B05} If $(G,A,B)$ is a graphic presentation of the isotropic system $S = (L,V)$ then $Q(G;x) = TM(S;x)$.
\end{theorem}

\section{Generalizations}

The closed form of the interlace polynomial in Theorem \ref{thm:vertexnullity} lends itself directly to a generalization to square matrices, and, by way of the adjacency delta-matroid of a graph, to delta-matroids. In each case, generalized pivot operations can be obtained that yield generalizations of the original recurrence for the interlace polynomial of a graph.

\subsection{Square matrices}

The adjacency matrix formulation of Aigner and van der Holst in Theorem \ref{thm:vertexnullity} lends itself nicely to a generalization of the polynomial to other matrices.

\begin{definition} \label{def:matrix-interlace} Let $A$ be a $V \times V$ matrix over the field $\mathbb{F}$. The \emph{interlace polynomial of $A$} is
    \begin{equation*}
        q_m(A;x) = \sum_{T \subseteq V} (x-1)^{n(A[T])}.
    \end{equation*}
\end{definition}

The recursive definition of the interlace polynomial for graphs can be recovered for general matrices using the following matrix operation, which has been extensively studied both in this context and others \cite{BH2011pivot,GP06,BH2011ni,T00}

\begin{definition} \label{def:ppt} Let $A$ be a $V \times V$ matrix over $\mathbb{F}$. Let $T \subseteq V$ such that the submatrix $A[T]$ is invertible over $\mathbb{F}$. There is a permutation matrix $X$ such that $XAX^T = \left( \begin{matrix} P & Q \\ R & S \end{matrix} \right)$ with $P$ the $T \times T$ submatrix of $A$. Then the principal pivot transform $A*T$ is the matrix satisfying \begin{equation*}X(A*T)X^T = \left( \begin{matrix} P^{-1} & -P^{-1}Q \\ RP^{-1} & S - RP^{-1}Q \end{matrix} \right).\end{equation*}
\end{definition}

The principal pivot transform can be thought of as a partial inverse.

\begin{theorem} \cite{T00} Let $A$ be an $n \times n$ matrix over a field $\mathbb{F}$ and let $T \subseteq \{1,\ldots,n\}$ such that $A[T]$ is invertible. Let $x$ and $y$ be vectors in $\mathbb{F}_n$. Let $u$ be the vector that agrees with $y$ on those entries indexed by $T$ and agrees with $x$ elsewhere. Let $v$ be the vector that agrees with $x$ on those entries indexed by $T$ and agrees with $y$ elsewhere. Then $A*T$ is the unique matrix satisfying
    \begin{equation*}
        y = Ax \text{ if and only if } (A*T)u = v
    \end{equation*}
for all vectors $x$ and $y$.
\end{theorem}

The first part of the following theorem relates the principal pivot transform to the pivot operation on a graph. The second part establishes that the interlace polynomial of a matrix satisfies a recurrence that, by the first part, generalizes the recurrence for the interlace polynomial of a graph. We will write $A \setminus v$ for the submatrix $A_{V \setminus \{v\}}$.

\begin{theorem} \label{thm:adjacencypivot} \
\begin{enumerate}
    \item \cite{BH2011pivot} Let $G$ be a graph with adjacency matrix $A$. Let $ab \in E(G)$. Then $G^{ab}$ has adjacency matrix $A * ab$ with the labels of $a$ and $b$ exchanged. 
    \item \cite{BH2011ni} Let $A$ be a $V \times V$ matrix over $\mathbb{F}$. Let $T \subseteq V$ with $A[T]$ invertible over $\mathbb{F}$. Then 
        \begin{equation*}q_m(A) = q_m(A\setminus v) + q_m((A*T) \setminus v)
        \end{equation*}
    for all $v \in V$.
\end{enumerate}
\end{theorem}

\subsection{Delta-matroids}

In \cite{BH2014ipm} Brijder and Hoogeboom generalized the interlace polynomial to delta-matroids, combinatorial objects that generalize matroids. In fact, they generalized the interlace polynomial to multimatroids, a further generalization of matroids introduced by Bouchet in a series of papers \cite{B97, B98, B01} which also generalize isotropic systems. We restrict here to the setting of delta-matroids, as this case most closely generalizes the interlace polynomial for graphs. Note, however, that even in this case the proofs of the theorems below often use the theory of multimatroids. We require first the following basic definitions regarding set systems.

\begin{definition} \emph{(set system)} \label{def:set-system} A \emph{set system} is a pair $(E,\mathcal{F})$ where $E$ is a finite set and $\mathcal{F} \subseteq 2^E$. The elements of $\mathcal{F}$ are called \emph{feasible sets}. A set system $(E,\mathcal{F})$ is said to be \emph{proper} if $E \neq \emptyset$. 
\end{definition}

\begin{definition} \emph{(delta-matroid)} \label{def:delta-matroid} A \emph{delta-matroid} is a proper set system $D = (E,\mathcal{F})$ satisfying the following \emph{symmetric exchange axiom}: for all $X,Y \in \mathcal{F}$, if $x \in X \Delta Y$ then there exists $y \in X \Delta Y$ such that $X \Delta \{x,y\} \in \mathcal{F}$.
\end{definition}

\begin{definition} \label{def:del-contract} Let $M = (E,\mathcal{F})$ be a set system. An element contained in every feasible set is a \emph{coloop}, and an element contained in no feasible set is a \emph{loop}. Let $e \in E$. If $e$ is not a coloop define \emph{$M$ delete $e$} to be the set system $M \setminus e = (E, \{F \in \mathcal{F}: e \not \in F\})$. If $e$ is not a loop define \emph{$M$ contract $e$} to be the set system $M/e = (E, \{F \setminus e: F \in \mathcal{F}, e \in F\})$. If $e$ is a coloop, define $M \setminus e = M/e$ and if $e$ is a loop define $M/e = M \setminus e$.
\end{definition}

\begin{definition} \label{def:twist} Let $M = (E,\mathcal{F})$ be a set system. For $X \subseteq M$, define the \emph{twist} $M*X$ to be the set system $(E, \{F \Delta X : F \in \mathcal{F}\})$.
\end{definition}

\begin{definition} \label{def:loopcomplement} Let $M = (E,\mathcal{F})$ be a set system. For $e \in E$, define the \emph{loop complement} $M+e$ to be the set system $(E, \mathcal{F} \Delta \{F \cup e: e \not \in \mathcal{F}\})$.
\end{definition}

\begin{definition} \label{def:dualpivot} Let $M = (E,\mathcal{F})$ be a set system. Twist $*e$ and local complementation $+e$ on a point $e \in E$ are involutions that generate a group isomorphic to $S_3$ \cite{BH2011pivot}. The third involution is $*e+e*e = +e*e+e$. It is called the \emph{dual pivot} and denoted $\bar{*} e$.
\end{definition}

\begin{definition} \label{def:distance} Let $M = (E,\mathcal{F})$ be a set system. For $X \subseteq E$, define the \emph{distance from $X$ to $M$} to be $d_M(X) = \min\{F \Delta X : F \in \mathcal{F}\}$.
\end{definition}

Note that loop complement on distinct points commutes \cite{BH2014ipm}, so we define $M+X$ for $X \subseteq E$ to be the set system obtained by performing loop complements at each of the points of $X$ in any order. Furthermore, note that while twist, deletion, and contraction are all operations on delta-matroids, loop complement is not.

\begin{example} \cite{CMNR14} \label{ex:loopcomplement} Let $M = (\{a,b,c\}, \{abc,ab,ac,bc,b,c,\emptyset\})$. Then $M$ is a delta-matroid, but $M+a = (\{a,b,c\}, \{a,b,c,bc,\emptyset\})$ is not.
\end{example}

\begin{definition} \label{def:vfsafe} \cite{BH2014ipm} We say a delta-matroid $M$ is \emph{vf-safe} if applying any sequence of twists and loop complements to $M$ yields a delta-matroid.
\end{definition}

The generalization of the interlace polynomial to delta-matroids follows from the vertex-nullity formula by first associating a delta-matroid to a graph via its adjacency matrix in such a way that the distance defined in Definition \ref{def:distance} above corresponds to the desired nullity of the graph.

\begin{definition} \label{def:adjacencydeltamatroid} Let $G$ be a graph with adjacency matrix $A$, considered over $GF(2)$. The \emph{adjacency delta-matroid of $G$}, denoted $M_G$, is the delta-matroid with ground set $V(G)$ and feasible sets consisting of all $X \subseteq V$ such that the principle submatrix $A[X]$ is invertible over $GF(2)$. Note that by convention $A[\emptyset]$ is invertible.
\end{definition}

We note that in the above definition, $GF(2)$ can be replaced with $GF(n)$. A delta-matroid $D=(E,\mathcal{F})$ is said to be \emph{representable over $GF(n)$} if for some $X \subseteq E$ there exists a skew-symmetric matrix $A$ over $GF(n)$ with $D = M_A * X$. This generalizes representability for matroids, and has been studied in \cite{B87delta, BD1991} among others. The following result shows that under this construction, distance for delta-matroids generalizes nullity for graphs.

\begin{theorem} \cite{BH2011ni} Let $G$ be a graph with adjacency delta-matroid $M_G$. Then $d_{M_G}(X) = n(G[X])$.
\end{theorem}

We can now define the interlace polynomial of a set system in such a way that, when the set system is the adjacency delta-matroid of a graph, it coincides with the interlace polynomial of a graph. Note that in \cite{BH2014ipm}, this definition is obtained via an evaluation of a generalized transition polynomial for multimatroids.

\begin{definition} \cite{BH2014ipm} \label{def:setsystem-interlace} Let $M = (E,\mathcal{F})$ be a set system. The \emph{interlace polynomial of $M$} is
    \begin{equation*}
        q_\Delta(M;x) = \sum_{X \subseteq E} x^{d_M(X)}.
    \end{equation*}
\end{definition}

\begin{theorem} \cite{BH2014ipm} \label{thm:deltamatroid-graph} Let $G$ be a graph with adjacency matroid $M_G$. Then
    \begin{equation*}
        q(G;x) = q_\Delta(M_G;x-1).
    \end{equation*}
\end{theorem}

The interlace polynomial of a delta-matroid also satisfies a recurrence generalizing that of the interlace polynomial of a graph.

\begin{theorem} \cite{BH2014ipm} Let $D = (E,\mathcal{F})$ be a delta-matroid. Let $e \in E$ be neither a loop nor a coloop. Then 
    \begin{equation*}
        q_\Delta(D;x) = q_\Delta(D \setminus e;x) + q_\Delta( D*e \setminus e;x).
    \end{equation*} 
If $\emptyset \in \mathcal{F}$, then for any $X \subseteq E$ and $e \in X$ we have 
    \begin{equation*}
        q_\Delta(D;x) = q_\Delta(D \setminus e;x) + q_\Delta( D*X \setminus e;x).
    \end{equation*}
If every element of $E$ is either a loop or a coloop, then $q_\Delta(D) = (y+1)^{|E|}$.
\end{theorem}

Since the empty matrix is by convention invertible over $GF(2)$, the adjacency delta-matroid of a graph always has $\emptyset$ feasible, and so the second recurrence above (which most directly generalizes the recurrence for graphs) holds. Moreover, the theorem below shows that this recurrence coincides precisely with the recurrence for graphs in the case that $D$ is the adjacency delta-matroid of a graph.

\begin{theorem} \cite{B87delta, BH2011pivot} Let $G$ be a graph with adjacency matrix $A$. Let $X \subseteq G$. If $A*X$ is defined, denote by $G*X$ the graph with adjacency matrix $A*X$. Then $M_{G*X} = M_G * X$ and $M_{G+X} = M_G + X$.
\end{theorem}

The $Q$ polynomial of a graph can also be generalized to delta-matroids.

\begin{definition} \cite{BH2014ipm} \label{def:globaldelta} Let $M = (E, \mathcal{F})$ be a set system. Define
    \begin{equation*}
        Q_\Delta(M;x) = \sum_{X \subseteq E} \sum_{Z \subseteq X} x^{d_{M+Z}(X)}.
    \end{equation*}
\end{definition}

\begin{theorem} \cite{BH2014ipm} \label{thm:globaldelta} Let $G$ be a simple graph. Then
    \begin{equation*}
        Q_\Delta(M_G;x-2) = Q(G;x).
    \end{equation*}
\end{theorem}

The recurrence for the polynomial $Q$ of graphs also generalizes. Here it is important to restrict to vf-safe delta-matroids, since loop complement is not an operation on general delta-matroids.

\begin{theorem} \cite{BH2011ni} Delta-matroids representable over $GF(2)$ (including adjacency delta-matroids) are vf-safe.
\end{theorem}

\begin{theorem} \cite{BH2014ipm} Let $D$ be a vf-safe delta-matroid. Then     \begin{equation*}
        Q_\Delta(D;x) = Q_\Delta(D \setminus e;x) + Q_\Delta(D *e \setminus e;x) + Q_\Delta(D \bar{*}e\setminus e;x)
    \end{equation*} for any $e \in E$ such that $e$ is neither a loop nor a coloop in $D$, and $e$ is not a coloop in $D\bar{*}e$.
\end{theorem}

The two-variable interlace polynomial of a graph can also be extended to delta-matroids as in the following definition and theorem.

\begin{definition} Let $M = (E,\mathcal{F})$ be a nonempty set system. Define 
    \begin{equation*}
        \bar{q}(M;x,y) = \sum_{X \subseteq E} x^{|X|} (y-1)^{n(X)}.
    \end{equation*}
\end{definition}

\begin{theorem} \cite{BH2014ipm} Let $G$ be a graph with adjacency delta-matroid $M_G$. Then
    \begin{equation*}
        \bar{q}\left(M_G;x-1,\frac{y-1}{x-1}\right) = q(G;x,y).
    \end{equation*}
\end{theorem}

The two-variable interlace polynomial of a delta-matroid also satisfies the following recurrence.

\begin{theorem} \cite{BH2014ipm} Let $D = (E,\mathcal{F})$ be a delta-matroid. Let $u \in E$. If $u$ is neither a loop nor a coloop, then
    \begin{equation*}
        \bar{q}(D;x,y) = \bar{q}(D \setminus u;x,y) + x \bar{q}(D*u\setminus u;x,y).
    \end{equation*}
If $u$ is a coloop, then
    \begin{equation*}
        \bar{q}(D;x,y) = (x+y) \bar{q}(D * u \setminus u;x,y),
    \end{equation*}
while if $u$ is a loop we have
    \begin{equation*}
        \bar{q}(D;x,y) = (1 + xy)\bar{q}(D \setminus u;x,y).
    \end{equation*}
\end{theorem}

Many evaluations of both the interlace polynomial and $Q$ for graphs extend to evaluations of the delta-matroid versions of these polynomials (and can often be obtained more easily in that context.) Note that item (6) below can be recovered from the Tutte-Martin polynomials \cite{B91}. Moreover, note that the evaluation of $Q(G)$ at $4$ does not extend to $Q_\Delta(D)$ \cite{BH2014ipm}.

\begin{theorem} \cite{BH2014ipm} Let $D = (E,\mathcal{F})$ be a delta-matroid with $n = |E|$. Then
\begin{enumerate}
    \item $q_\Delta(D;1) = 2^n$;
    \item  $q_\Delta(D;0) = |\mathcal{F}|$;
    \item  if all sizes of feasible sets in $D$ have the same parity,    then $q_\Delta(D;-1) = 0$;
    \item  if $D$ is vf-safe then $Q_\Delta(D;-2) = 0$;
    \item  if $D$ is vf-safe then $q_\Delta(M;-2) = (-1)^n(-2)^{d_{D\bar{*}E}}$; and
    \item if $D$ is binary (i.e. representable over $GF(2)$) then $q_\Delta(D;2) = kq_\Delta(D;-2)$ for some odd integer $k$.
\end{enumerate}
\end{theorem}

Lastly, we remark that there is a connection, as with the interlace polynomial of graphs, between the interlace polynomial of delta-matroids restricted to matroids and the Tutte polynomial along $x=y$.

\begin{definition} \cite{BH2014ipm} Let $M = (E,\mathcal{F})$ be a matroid described by its bases (i.e. a delta-matroid where all feasible sets have the same cardinality). The Tutte polynomial $t(M;x,y)$ of $M$ is given by the recurrence
\begin{enumerate}
    \item $t(M;x,y) = t(M/e;x,y) + t(M \setminus e;x,y)$ if $e$ is neither a loop nor a coloop of $M$, and
    \item $t(M;x,y) = x^i y^j$ if $M$ consists of $i$ coloops and $j$ loops.
\end{enumerate}
\end{definition}

\begin{theorem} \label{thm:tuttematroid}\cite{BH2014ipm} Let $M = (E,\mathcal{F})$ be a matroid. Then 
    \begin{equation*}
        t(M;x,x) = q_\Delta(M; x-1).
    \end{equation*}
\end{theorem}

\section{Conclusion}

We collect in Tables \ref{tab:graph} and \ref{tab:4regular} the known connections between polynomials and combinatorial objects considered above.

There are a number of remaining research directions regarding the interlace polynomials. A natural question is whether results on the interlace polynomial of $4$-regular graphs and two-in two-out digraphs can be extended to arbitrary Eulerian graphs and digraphs. Transition systems and the Martin polynomials both extend to this case (see e.g. \cite{EM98}), but no results are known on the interlace polynomial. 

There are also further directions of research involving the Tutte polynomial. Brijder and Hoogeboom found in \cite{BH2014ipm} that a generalization of the two-variable interlace polynomial for matroids captures more of the Tutte polynomial than is captured in Theorem \ref{thm:tuttematroid}, raising the question of whether there is some general combinatorial object and variant of the interlace polynomial that capture the entirety of the Tutte polynomial. Note that the recursive relation defining the Tutte polynomial is not well-defined on delta-matroids, so it is likely that another context is needed.

In \cite{BR01}, Bollob\'{a}s and Riordan generalized the Tutte polynomial to a polynomial of embedded graphs, which has been shown in \cite{EMS07} to have a connection to the two-variable interlace polynomial analogous to the connection of Theorem \ref{thm:tutteplane}. This demonstrates that the interlace polynomial contains some topological information, and a natural question to ask is whether there is a full extension of the interlace polynomial to embedded graphs.

Lastly, in the case of the interlace polynomials of graphs, the study of $q_N$ and $Q$ has dominated the literature, and the properties and structure of the two-variable polynomial and the vertex-rank polynomial are less well-known. Furthermore, the closed forms of the vertex-rank and vertex-nullity polynomials suggest the possibility of defining related polynomials with respect to the incidence or Laplacian matrices of a graph, and studying the general theory of such vertex-rank/nullity polynomials.

\begin{table}[h]
\centering
\begin{tabular}{|c||c|c|c|}
\textbf{Graph} & \multicolumn{3}{|c|}{\textbf{Combinatorial objects}} \\ \cline{2-4}
  \textbf{polynomials} & \textbf{Isotropic system} & \textbf{Matrix} & \textbf{Delta-matroid} \\ \hline \hline
 A graph $G$ & $(G,A,B)$ & $A(G)$ & $M_G$ \\ \hline
 $q(G;x)$ & $tm(S,A+B;x)$ & $q_m(A(G);x)$ & $q_\Delta(M_G;x-1)$ \\ \hline
 $q(G;x,y)$ & & & $\bar{q}\left(M_G; x-1, \frac{y-1}{x-1} \right)$ \\ \hline
 $Q(G;x)$ & $TM(S;x)$ & & $Q_\Delta(M_G;x-2)$
\end{tabular}
\caption{The first column gives the various interlace polynomials of a graph $G$. Reading across gives the combinatorial objects generalizing graphs, how graphs are encoded by each, and the specializations of polynomials of more general objects that capture the interlace polynomial of a graph in question.}
\label{tab:graph}
\end{table}

\begin{table}[h]
\centering
\begin{tabular}{|c||c|c|c|}
\textbf{Graph}& \multicolumn{3}{|c|}{\textbf{Types of graph $G$.}} \\ \cline{2-4}
\textbf{polynomials} & \textbf{$4$-regular} & \textbf{Two-in two-out} & \textbf{Plane} \\ \hline \hline
$m(G;x)$ & & $q(H(C);x)$ &  \\ \hline
$M(G;x)$ & $Q(H(C);x)$ & &  \\ \hline
\vbox{\hbox{\strut $t(G;x,x)$}\hbox{\strut }} & & & \vbox{\hbox{\strut $q(H(C);x)$}\hbox{\strut $m(\vec{G}_m;x)$}} \\
\end{tabular}
\caption{The first column gives specializations of the Martin and Tutte polynomials of a graph $G$. Reading across gives the appropriate form of interlace polynomial capturing those specializations in the case that $G$ is $4$-regular, two-in two-out, or plane. In the table, $C$ is an Eulerian circuit of $G$.}
\label{tab:4regular}
\end{table}

\bibliographystyle{plain}
\bibliography{interlace_bibtex}

\end{document}